\documentclass{article}

\newtheorem{thm}{Theorem}[section]
\newtheorem{Theorem}{Main Result}
\newtheorem{rem}[thm]{Remark}
\newtheorem{prop}[thm]{Proposition}
\newtheorem{lemma}[thm]{Lemma}

\def\qed{{\hfill\hphantom{.}\nobreak\hfill$\Box$}}
\usepackage{amssymb,amsmath,hyperref}
\usepackage[english]{babel}

\newcommand{\dd}{\mathsf{d}}
\newcommand{\R}{\mathbb{R}}
\newcommand{\N}{\mathbb{N}}
\newcommand{\Isom}{\mathrm{Isom}}
\newcommand{\Res}{\mathrm{Res}}
\newcommand{\proof}{\emph{Proof.~}}

\begin{document}

\title{Rigidity at infinity of trees and Euclidean buildings}

\author{Koen Struyve\thanks{The author is supported by  the Fund for Scientific Research -- Flanders (FWO - Vlaanderen).}
}



\maketitle

\begin{abstract}
We show that if a group $G$ acts isometrically on a locally finite leafless $\mathbb{R}$-tree inducing a two-transitive action on its ends, then this tree is determined by the action of $G$ on the boundary. As a corollary we obtain that locally finite irreducible Euclidean buildings of dimension at least two are determined by their complete building at infinity.
\end{abstract}


\section{Introduction}
Euclidean buildings (also known as $\mathbb{R}$-buildings or affine apartment systems) form, together with symmetric spaces, an important class of CAT(0)-spaces. The boundary at infinity of a Euclidean building comes equipped with a simplicial structure making it into a spherical building. In the case of one-dimensional Euclidean buildings, which are exactly the leafless $\mathbb{R}$-trees, this boundary is just a set.

We prove two rigidity results dealing with the question how this structure at infinity determines a tree or a Euclidean building. The first result is  for leafless locally finite $\mathbb{R}$-trees with a two-transitive action at infinity.

\begin{Theorem}\label{thm:1}
Let $G$ be a group acting isometrically on two  leafless trees $T_1$ and $T_2$, and assume that there is a $G$-equivariant bijective map $f : \partial_\infty T_1 \to \partial_\infty T_2$. 

If $T_1$ and $T_2$ are  both locally finite, admit at least three ends and $G$ acts two-transitively on $\partial_\infty T_1$, then there exists (after rescaling the metric on $T_2$) a isometry $\bar{f}: T_1 \to T_2$ which induces the map $f$ on $\partial_\infty T_1$. 

This isometry is unique if $T_1$ has at least two branch points.
\end{Theorem}  

Using this we then prove the following rigidity result for Euclidean buildings. 

\begin{Theorem}\label{thm:2}
If $X_1$ and $X_2$ are two locally-finite irreducible Euclidean buildings of dimension at least two with an isometry $f: \partial_\infty X_1 \to \partial_\infty X_2$ between their buildings at infinity, then there exists (after rescaling the metric on $X_2$) an isometry $\bar{f}: X_1 \to X_2$ which induces the map $f$ on $\partial_\infty X_1$. This isometry is unique if $X_1$ is not an Euclidean cone.
\end{Theorem}

These Main Results extend the work of B. Leeb (see~\cite[Prop. 4.20 \& Addendum  1.3]{Lee:00}), who obtained a similar result under the additional assumption that $f$ is continuous for the cone topology, and the work of L. Kramer and R. Weiss, which prove a rigidity result for trees and Euclidean buildings concerning quasi-isometries (\cite{Kra-Wei:*}, see also~\cite{Kle-Lee:97}).

In the case of Bruhat-Tits buildings (i.e. Euclidean buildings of algebraic origin) such rigidity results were already known for algebraic reasons (see~\cite[27.6]{Wei:09}), based on a result of F. K. Schmidt that a field with multiple complete valuations on it is necessarily algebraically closed (see~\cite{Sch:33}), which implies a non-discrete value group and an infinite residue field.

In~\cite{Sch:33} one also obtains the existence of fields with multiple complete valuations, implying the existence of counterexamples in the non-discrete case if we drop the local finiteness condition.

Finally we remark that Main Result~\ref{thm:1} does not hold if one would replace local finiteness by discreteness, as one can construct a counterexample using transfinite recursion (see~\cite{Str:*}). For Main Result~\ref{thm:2} no such counterexamples are known at the moment.

\section{Note about the proof}

The main innovation of this paper is the recognition of the bounded subgroups by the action on the boundary at infinity of the tree (see Proposition~\ref{prop:rec}). After this point we follow the same arguments as used in the proof of L. Kramer and R. Weiss (\cite{Kra-Wei:*}). (Alternatively one could follow the arguments in~\cite{Lee:00}.)


\section{Trees and actions} \label{section:trees}

In this section we collect some basic results on trees and isometries acting on these. For a detailed discussion we refer to~\cite{Chi:01} and~\cite{Mor-Sha:84}.

\subsection{Definitions}
A metric space $(T,\dd)$ is an \emph{$\R$-tree}, or shortly a tree, if it satisfies the following two properties.
\begin{itemize}
\item[(T1)] $(T,\dd)$ is a uniquely geodesic metric space (any two points are joined by a unique geodesic segment).
\item[(T2)] If two geodesic segments meet only at a common end-point, then their union is a geodesic segment.
\end{itemize}

We say that $(T,\dd)$ is \emph{leafless} if it has extensible geodesics. A \emph{ray} of a tree is an embedded closed half-line. Two rays are \emph{equivalent} if they intersect in a ray. The equivalence classes are called the \emph{ends} of $T$, the set of all ends is denoted by $\partial_\infty T$. 

If $a$ and $b$ are two different ends in $\partial_\infty T$, then we define $(a,b)$ to be the unique isometric image of $\R$ in $T$ containing rays with ends $a$ and $b$. We call this the \emph{apartment} with ends $a$ and $b$. 

Two rays starting in the same point $x\in T$ \emph{locally coincide} if they intersect in more than just $x$. This forms an equivalence relation between rays, the set of equivalence classes is called the \emph{space of directions at $x$}. The \emph{valency} of a point in $T$ is the cardinality of its space of directions.   

A \emph{branch point} is a point of valency at least three. If there are no points with infinite valency, we say that $T$ is \emph{locally finite}. Note that locally finite does not necessarily imply locally compact.

\subsection{Isometries of trees}

Let $\Isom(T)$ be the isometry group of the leafless tree $T$ (acting from the left). This group has an induced action on $\partial_\infty T$.

For $g \in \Isom(T)$ define the \emph{minimal displacement length} $l(g)$ as $\inf_{x \in T} \dd(x,gx)$. We denote by $A_g$ the set $\{x\in T \vert \dd(x,gx) = l(g)\}$ (see~\cite[p. 82]{Chi:01}). One easily verifies that $A_{g^{-1}} = A_g$ and $A_{hgh^{-1}} = h A_g$, with $h \in \Isom(T)$.

An isometry $g$ of $T$ is either \emph{elliptic}, in which case $g$ fixes a point in $T$ (so $A_g$ consists of the fixed points of $g$), or it is \emph{hyperbolic}, in which case $A_g$ is an apartment such that the restriction of $g$ to $A_g$ is a translation of length $l(g)$ and such that the fixed ends of $g$ in $\partial_\infty T$ are exactly the ends of $A_g$.

Note that for an hyperbolic isometry $g$ one has for $n \in \mathbb{Z} \setminus \{0\}$ that $A_g = A_{g^n}$. 

\subsection{Pairs of isometries}

In this section we study how two isometries of a tree interact.

\begin{lemma}\label{lemma:2iso}
If $g,h \in \Isom(T)$ are two isometries such that the intersection of $A_g$ and $A_h$ is empty, then $gh$ is hyperbolic and $A_{gh}$ intersects both $A_g$ and $A_h$. Moreover $l(gh) = l(g)+l(h) +2\dd(A_g, A_h) $. 
\end{lemma}
\proof
This follows directly from~\cite[Lem. III.2.2]{Chi:01}. \qed

\begin{lemma}\label{lemma:test}
If $g,h \in \Isom(T)$ are two isometries such that $h$ is hyperbolic while $g$ and $gh^n$ are elliptic for every $n \in \mathbb{Z}$, then $g$ stabilizes $A_h$.
\end{lemma}
\proof
Fix $n$ to be an even number $2m$ different from zero. The sets $A_{g}$ and $A_h$ intersect by Lemma~\ref{lemma:2iso} as $gh$ is elliptic. By the same lemma the sets $A_{gh^n}$ and $A_{h^{-n}}= A_h$ intersect as $g$ is elliptic. Let $x$ be a point in the intersection of $A_{gh^n}$ and $A_h$. 

The isometry $gh^n$ translates the apartment $A_h$ by $\vert n \vert l(h)$ and then applies $g$ to it. As $g$ fixes some point of $A_h$, the only possibility is that $g$ fixes a unique point $y$ of $A_h$ which is the midpoint of the geodesic segment $[x, h^n x]$. In particular we have that $y =h^m x$ and that $g$ maps $[h^{-m}y, y]$ to $[y, h^m y]$. 

As this holds for arbitrary $m \in \mathbb{Z} \setminus \{ 0 \}$ it follows that $g$ stabilizes $A_h$.
\qed

\section{Euclidean buildings}

In this section we discuss $\R$-buildings. For a detailed treatment we refer to~\cite{Par:00} and~\cite{Tit:86}. 

\subsection{Definitions}\label{section:EuclDef}
Let $(\overline{W},S)$ be a spherical Coxeter system of rank $n$ (with $n \geq 1$). The group $\overline{W}$ has a standard representation as a finite reflection group acting on a $n$-dimensional Euclidean space $\mathbb{A}$, called the \emph{model space}. Let $W$ be the group acting on $\mathbb{A}$ generated by $\overline{W}$ and the translations of $\mathbb{A}$.

The \emph{walls} of $\mathbb{A}$ are those hyperplanes fixed by conjugates in $W$ of involutions in $S$. A \emph{root} is a (closed) half-space of $\mathbb{A}$ bordered by a wall. The set of all walls of $\mathbb{A}$ through a given point defines a poset of simplicial cones in $\mathbb{A}$ (called \emph{sector-faces}), which forms the simplicial complex of the Coxeter system $(\overline{W},S)$. The maximal cones are called \emph{sectors}. (See~\cite[Chapter 1]{Abr-Bro:08} for a detailed discussion on finite reflection groups.)

Let $(X,\dd)$ be a metric space together with a collection $\mathcal{F}$ of isometric injections (called \emph{charts}) from the model space $\mathbb{A}$ into $X$. An image of the model space is called an \emph{apartment}, an image of a root a \emph{half-apartment} and an image of a sector(-face) is called again a \emph{sector(-face)}. The space $X$ together with the collection $\mathcal{F}$ is a \emph{Euclidean building} if the following 5 properties are satisfied:
\begin{itemize}
 \item[(A1)] If $w\in W$ and $f\in \mathcal{F}$, then $f \circ w \in \mathcal{F}$.
 \item[(A2)] If $f,f' \in \mathcal{F}$, then $X :=f^{-1}(f'(\mathbb{A}))$ is a closed and convex subset of $\mathbb{A}$, and $f|_X = f'\circ w|_X$ for some $w\in W$.
 \item[(A3)] Each two points of $X$ lie in a common apartment.
 \item[(A4)] Any two sectors $S_1$ and $S_2$ contain subsectors $S_1' \subset S_1$ and $S_2' \subset S_2$ lying in a common apartment.
 \item[(A5)] If three apartments intersect pairwise in half-apartments, then the intersection of all three is non-empty.
\end{itemize}

The \emph{dimension} of a Euclidean building is the dimension of the model space $\mathbb{A}$. One-dimensional Euclidean buildings are exactly the leafless $\R$-trees.

We say that a Euclidean building is \emph{irreducible} if is not a Euclidean space nor decomposes as a product of two Euclidean buildings (where the product is defined as in~\cite[I.5.1]{Bri-Hae:99}).

\subsection{Global and local structure}
Two sector-faces are \emph{parallel} if the Hausdorff distance between both is finite. This relation is an equivalence relation due to the triangle inequality. The equivalence classes, named \emph{simplices at infinity},  form  a spherical building of type $(\overline{W},S)$ called the \emph{building at infinity} of the Euclidean building. 

Different choices of the collection of charts lead to different buildings at infinity. There is however a unique maximal choice for the collection of charts, from which one obtains the \emph{complete building at infinity} $\partial_\infty X$. Throughout this paper we will always assume that we are in this situation. The metric realization of the complete building at infinity is isometric to the \emph{Tits boundary}  of $X$ as defined in~\cite[Chapter II.9]{Bri-Hae:99}.

One can also define local equivalences. Let $x$ be a point of $X$, and $F,F'$ two sector-faces based at $x$. Then these two sector-faces will \emph{locally coincide} if their intersection is a neighborhood of $x$ in both $F$ and $F'$. This relation forms an equivalence relation defining \emph{germs of faces} as equivalence classes. These germs form a (weak) spherical building of type $(\overline{W},S)$,
called the \emph{residue} at $x$. \emph{Weak} means here that some panels might be contained in only two chambers.  

We say that an Euclidean building is \emph{locally finite} if every residue is finite as simplicial complex. Again we note that locally finite does not imply locally compact.

\subsection{Panel trees and projectivities}\label{section:paneltrees}

To every panel $a$ of $\partial_\infty X$, one can associate a \emph{panel-tree} (see~\cite[Prop. 4]{Tit:86})  denoted by $X(a)$. The ends of this tree are in bijective correspondence with the vertices of the residue $\Res(a)$ of the panel. 

If two panels $a$ and $b$ are opposite, then there exists a canonical isometry between the corresponding panel-trees, which at infinity (after identification with $\Res(a)$ and $\Res(b)$) induces the projection map between $\Res (a)$ and $\Res (b)$. This projection map is also called the \emph{perspectivity} between $\Res (a)$ and $\Res (b)$.

The isometry between the panel trees is obtained by constructing a wall tree associated to the unique wall of $\partial_\infty X$ containing both $a$ and $b$, and showing that this wall tree is isometric (in a canonical way) to both the panel trees $X(a)$ and $X(b)$  (see~\cite[\S 9]{Tit:86} for details).

One can chain any number of perspectivities together to obtain maps from $\Res(a)$ to itself, called \emph{projectivities}. On the level of the panel tree this yields isometries from $X(a)$ to itself. The set of all projectivities of $\Res(a)$ forms a group, called the \emph{projectivity group} of $a$. If $X$ is an irreducible Euclidean building of dimension at least two with a thick spherical building at infinity, then this projectivity group acts two-transitively on $\Res(a)$ (see~\cite[1.2]{Kna:88} and~\cite[Thm. 3.14]{Kra-Wei:*}). Hence, on the level of the panel tree, one obtains that the projectivity group acts as a group of isometries of this tree inducing a two-transitive action on its set of ends. (See~\cite[Section 3]{Lee:00} for a related concept.)

\subsection{Euclidean cones}\label{section:cones}

Let $(B,\dd_B)$ be the metric realization of a given spherical building and $X$ the quotient of $B \times [0,\infty[$ where one identifies the subset $B \times \{0\}$  to a point. The space $X$, with the metric $\dd((x,s), (y,t))^2 = s^2 + t^2 - 2st \cos(\dd_B (x,y))$ (see~\cite[I.5.6]{Bri-Hae:99}), is called the  \emph{Euclidean cone} over $B$, which in itself is a Euclidean building.

The number of branch points of a panel tree of such an Euclidean building is either 0 or 1, where the first possibility occurs if the  spherical building one started with is not thick.

\subsection{Thick reductions}\label{section:reductions}
Let $X$ be a Euclidean building. Following~\cite[Section 4.9]{Kle-Lee:97} and~\cite[4.12]{Kra-Wei:*}, one can reduce the spherical Coxeter group $W$ used in Section~\ref{section:EuclDef} to define the Euclidean building  to a minimal one, making its spherical building at infinity the product of a sphere $\mathbb{S}^k$ and a thick spherical building. This thick spherical building depends only on the original spherical building at infinity, not of the Euclidean building itself.

A consequence of this operation is that if an irreducible Euclidean building is not an Euclidean cone, then every panel tree of its reduction contains at least two branch points, see~\cite[4.26]{Kra-Wei:*}.

\section{Proof of Main Result~\ref{thm:1}}

In this section, assume that $T$ is a locally finite leafless tree, admits at least three ends and that we have a subgroup $G$ of $\Isom(T)$ inducing a two-transitive action on $\partial_\infty T$. We will aim to reconstruct the tree from $\partial_\infty T$ and the action of $G$.

We first mention two results from~\cite{Kra-Wei:*} concerning trees with a two-transitive action at infinity .

\begin{lemma}\label{lemma:types}
The set of branch points of $T$ satisfies one of the following possibilities.
\begin{itemize}
\item[Type (I)] There is a single branch point, and $T$ is the Euclidean cone over its set of ends $\partial_\infty T$.
\item[Type (II)] There is an infinite discrete set of branch points, and $T$ is a simplicial metric tree where every vertex has valency at least 3, and all edges have the same length.
\item[Type(III)] The set of branch points is dense, in particular the set of branch points in a single apartment is dense.
\end{itemize} 
\end{lemma}
\proof
This is essentially~\cite[Prop 2.5 \& Cor. 2.6]{Kra-Wei:*}.
\qed

\begin{lemma} \label{lemma:2trans}
Let $x$ be a branch point of $T$. The stabilizer $G_x$ of this point induces a two-transitive action on the space of directions at $x$.  
\end{lemma}
\proof
This follows directly from~\cite[Prop. 2.3]{Kra-Wei:*}.
\qed

By Lemma~\ref{lemma:types} and the local finiteness of $T$, we observe that Type (I) occurs if and only if the set of ends $\partial_\infty T$ is finite. In this case it is trivial to reconstruct $T$, so we can assume that $T$ is either of Type (II) or (III). 

\begin{lemma}\label{lemma:spider}
Let $x$ be a branch point of $T$, then there exists an isometry in $G$ fixing only $x$. In particular this isometry acts freely on $\partial_\infty T$.
\end{lemma}
\proof
The group $G_x$ acts two-transitive on the space of directions at $x$. Let $\overline{G_x}$ be the normal subgroup of $G_x$ consisting of those isometries acting trivially on this space of directions. 

By the orbit-counting theorem (see for example~\cite{Neu:79}) the average number of fixed directions of the elements in the finite group $G_x / \overline{G_x}$ is 1. As the identity fixes at least three directions, there exists an element in $G_x$ fixing no directions at $x$. Such an isometry clearly fixes only $x$ and cannot fix any end.
\qed

The next proposition allows us to recognize the elliptic elements in $G$.

\begin{prop}\label{prop:ell}
An isometry $h \in G$ is elliptic if and only if one of the following two conditions is satisfied:
\begin{itemize}
\item
The order of the set of fixed ends by $h$ in $\partial_\infty T$ is different from two.
\item
$h$ fixes exactly two elements $a$ and $b$ of $\partial_\infty T$ and  there exists an isometry $g \in G$ acting freely on $\partial_\infty T$, not stabilizing $\{a,b\}$ and an $N \in \mathbb{N} \setminus \{0\}$ such that for every $n \in \mathbb{Z} \setminus \{0\}$ the isometry $gh^{Nn}$ acts freely on $\partial_\infty T$.
\end{itemize}
\end{prop}
\proof 
We first proof the ``if'' part. Assume by way of contradiction that $h$ is hyperbolic.  A hyperbolic isometry fixes exactly two ends (which we set to be $a$ and $b$), so we can assume that the second condition is satisfied for a certain $g$ and $N$. 

The free action on $\partial_\infty T$ implies that the isometries $g$ and $gh^{Nn}$ ($n \in \mathbb{Z} \setminus \{0\}$) are elliptic. As $h^N$ is still hyperbolic, Lemma~\ref{lemma:test} yields that $g$ stabilizes $A_{h^N}$, so it stabilizes its set of ends $\{a,b\}$. Hence we have obtained a contradiction.

We now prove the ``only if'' part. Assume that $h$ is an elliptic isometry in $G$. Without loss of generality we may assume that $h$ fixes exactly two elements $a$ and $b$ of $\partial_\infty T$. This implies that $A_h$ contains the apartment $(a,b)$. Let $x$ be a branch point on this apartment and $m$ the number of directions at $x$. Set $N$ to be $m$ factorial. The power $h^N$, as well as all its subsequent powers, fixes all directions at $x$ as $N$ is the order of the symmetric group on $m$ symbols. 

Let $g$ be an isometry in $G$ fixing only $x$, as constructed in Lemma~\ref{lemma:spider}. One can pick $g$ in such a way that it does not interchange the ends $a$ and $b$, for example by conjugating $g$ with an element in $G_x$ fixing $a$ but not $b$. 

The product $gh^{Nn}$ ($n \in \mathbb{Z} \setminus \{0\}$) fixes $x$, but none of the directions at $x$, hence it fixes exactly $\{x\}$ and acts freely on $\partial_\infty T$. We conclude that $h$ satisfies the second condition. \qed

A subgroup $H$ of $G$ is \emph{bounded} if at least one orbit, and hence every orbit, of $H$ on points of $T$ has bounded diameter.

\begin{prop}\label{prop:rec}
One can recover the bounded subgroups of $G$ from the action of $G$ on $\partial_\infty T$.
\end{prop}
\proof
We observe that bounded subgroups consist of elliptic elements only. 

Suppose that $H < G$ is a subgroup which consists of elliptic elements only, we then have to determine if $H$ is bounded, or equivalently if $H$ fixes a point of the metric completion $\overline{T}$ (by the Bruhat-Tits fixed point theorem~\cite[Prop. 3.2.4]{Bru-Tit:72}). 

If $H$ is finitely generated  then it fixes an element of $T$ by~\cite[Prop.~II.2.15]{Mor-Sha:84}. For the general case consider the fixed sets of the finitely generated subgroups of $H$. These form a filtering family of closed convex subsets of  $\overline{T}$. Such a filtering family has a common element in $\overline{T}$ or $H$ fixes a subset of intrinsic radius at most $\pi/2$ in its boundary by~\cite[Thm.~1.1]{Cap-Lyt:09}. The first possibility implies that $H$ is bounded, the second that $H$ fixes exactly one end of $T$.

Now suppose the second case holds. If $H$ would fix a point $x$  of $\overline{T}$, then it would also fix a branch point $y$ in $T$ by considering a ray with base point $x$ in the equivalence class of the fixed end of $H$. (Such a branch point $y$ exists by Lemma~\ref{lemma:types}, and the assumption that $T$ is of Type (II) or (III).) The stabilizer of $y$ contains $H$, does not admit fixed ends (by Lemma~\ref{lemma:spider}), and consists of elliptic elements only, hence we clearly can recognize it as being bounded. 

As subgroups of bounded subgroups are again bounded, we can hence recognize the bounded subgroups of $G$ as those subgroups contained in a subgroup of $G$ consisting of elliptic elements only and fixing no ends of $T$. \qed

\begin{rem}
Alternatively one can try to generalize~\cite[Prop. 3.4]{Tit:70}.
\end{rem}

Proposition~\ref{prop:rec} and~\cite[Lem. A.4]{Kra-Wei:*} allows us to recover the maximal bounded subgroups corresponding with the $G$-isolated points of $T$. Using~\cite[(2.10)-(2.13)]{Kra-Wei:*} and its appendix, one is then able to reconstruct the structure of the tree $T$ up to a scaling factor. Hence Main Result~\ref{thm:1} follows.


\section{Proof of Main Result~\ref{thm:2}}
Let $X_1$, $X_2$ and $f$ be as in the statement of Main Result~\ref{thm:2} and assume that the buildings at infinity of $X_1$ and $X_2$ are thick. If this is not the case we can reduce it to the thick case, as in Section~\ref{section:reductions}.

According to Section~\ref{section:paneltrees} the projectivity group associated to a panel $a$ of $\partial_\infty X_1$ acts isometrically on $X_1(a)$ while inducing a two-transitive action on its set of ends. So we may apply Main Result~\ref{thm:1} and obtain (after rescaling) an isometry $\bar{f}_a$ from $X_1(a)$ to $X_2(f(a))$ inducing $f\vert_{\Res(a)}$ on its set of ends. 

We now want to apply~\cite[Thm. 2]{Tit:86} (see also~\cite[Thm. 6.17]{Hit:09} and~\cite[Section 5.5]{Lee:00}) which would yield that there exists, after rescaling the metric on $X_2$, an isometry $\bar{f} :X_1 \to X_2$ which induces the map $f$ on $\partial_\infty X_1$.

However one has to take into account the possibility that the scaling factor depends on the choice of panel $a$. If two panels $a$ and $b$ of $\partial_\infty X_1$ can be connected via a series of perspectivities, one obtains via the methods of Section~\ref{section:paneltrees} that the scaling factors for the panel-trees $X_1(a)$ and $X_1(b)$ are the same.

This suffices in the proof of~\cite[Thm. 6.17]{Hit:09} to construct a bijection $\bar{f}$ from $X_1$ to $X_2$, mapping apartments to apartments and walls to walls. As the Euclidean building $X_2$ is irreducible, it follows that $\bar{f}$ is an isometry up to rescaling the metric on $X_2$. (To see this observe that the walls in an apartment define the metric in the apartment uniquely up to scaling.) Moreover the scaling factor will be unique if $X_1$ is not an Euclidean cone.

This concludes the proof of Main Result~\ref{thm:2}.

\appendix
\section{Addendum}
The goal of this addendum is to prove the following existence result showing that Main Result 1 fails when one does not require local finiteness, even in the discrete case. 

\begin{thm}\label{thm:3}
There exists a group $G$ acting isometrically on two  leafless discrete trees $T_1$ and $T_2$, such that there exists $G$-equivariant bijective map $f : \partial_\infty T_1 \to \partial_\infty T_2$ which does not extend (even after rescaling the metric on $T_2$) to an isometry $\bar{f}: T_1 \to T_2$.
\end{thm}  

The proof is constructive using techniques from transfinite recursion. While quite technical it is somewhat straightforward.

\subsection{Sequences and gluing trees}\label{section:treesq}

In this section we provide the machinery for the intuitive process of gluing trees together along isometric subtrees and working with sequences of trees. For this we use the theory of CAT($\kappa$)-spaces as described in~\cite{Bri-Hae:99}. The next lemma provides an equivalent definition of $\R$-trees.

\begin{lemma}\label{lemma:cat}
A geodesic metric space is a tree if and only if it is a CAT($\kappa$)-space for all $\kappa \in \R$.
\end{lemma}
\proof
See~\cite[Ex. 1.15.(5)]{Bri-Hae:99}.
\qed

Let $(T_1,\dd_1)$, $(T_2,\dd_2)$ and $(A,\dd)$ be three ($\R$-)trees together with isometric embeddings $\iota_1: A \to T_1$ and $\iota_2: A \to T_2$. By~\cite[Def. I.5.23]{Bri-Hae:99} one can consider the \emph{gluing $T_1 \sqcup_ A T_2$ of $T_1$ to $T_2$ along $A$}. This is a geodesic metric space admitting natural isometric embeddings $p_1: T_1 \to T_1 \sqcup_ A T_2$ and $p_2 : T_2 \to T_1 \sqcup_ A T_2$
 by~\cite[Lem. I.5.24]{Bri-Hae:99}. The gluing is again an $\R$-tree by~\cite[Thm. II.11.3]{Bri-Hae:99} and Lemma~\ref{lemma:cat}.

The next lemma deals with sequences of trees.

\begin{lemma}\label{lemma:treesq}
Let $I$ be a totally ordered set, with for each $i \in I$ an $\R$-tree $(T_i,\dd_i)$ such that if $i,j \in I$ and $i < j$, then $(T_i,\dd_i)$ is a metric subspace of $(T_j,\dd_j)$ (i.e. $T_i \subset T_j$ and $d_i = d_j\vert_{T_i}$). Consider the set $T = \cup_{i\in I} T_i$ with the metric $\dd$ induced by the metrics $\dd_i$ on the subspaces $T_i$ ($i \in I$). Then $(T,\dd)$ is again an $\R$-tree.
\end{lemma}
\proof
As each of the trees $(T_i,\dd_i)$ is a geodesic metric space, the tree $(T,\dd)$ is as well. In order to check that $(T,\dd)$ is a  CAT($\kappa$)-space for a given $\kappa$ one only needs to consider configurations of four points in $T$  (see~\cite[Prop. II.1.11]{Bri-Hae:99}). As one can always find an $i\in I$ such that $T_i$ contains a given 4-tuple of points in $T$, it follows that $(T,\dd)$ is a CAT($\kappa$)-space if each of the spaces $(T_i,\dd_i)$ ($i \in I$) is. Applying Lemma~\ref{lemma:cat} concludes the proof. \qed

An important consequence is that one can transfinitely recursively  glue trees together for any well-ordered set (see~\cite[I.5.26]{Bri-Hae:99} for countable sequences). Let $I$ be a well-ordered set with for each $i\in I$ a tree $(T_i,\dd_i)$. For each $i \in I$ we recursively define a tree $(\widetilde{T}_i, \tilde{\dd}_i)$ in the following way.

\begin{itemize}
\item If $i$ is an isolated point (i.e. it has a predecessor $j$), then we define $(\widetilde{T}_i, \tilde{\dd}_i)$ by gluing $(\widetilde{T}_j, \tilde{\dd}_j)$ and $(T_i, \dd_i)$ along certain isometric subtrees which are to be made precise in the actual construction.
\item If $i$ is a limit (i.e. it has no predecessors), then we obtain $(\widetilde{T}_i, \tilde{\dd}_i)$ by gluing to $(T_i, \dd_i)$ to the union of all $(\widetilde{T}_j, \tilde{\dd}_j)$ where $j\in I, j<i$ (which is a tree by Lemma~\ref{lemma:treesq}) in some way to be made precise in the actual construction.
\end{itemize}

The $\R$-tree resulting from this transfinite recursion is the union of all these $(\widetilde{T}_i, \tilde{\dd}_i)$.

\subsection{A construction}\label{section:example}

In this section we construct the group and the trees with the properties stated in Theorem~\ref{thm:3}.

We will assume that every tree in this section is discrete with edge length 1.


We say that a 4-tuple $(T,\dd,D,G)$, where $(T,\dd)$ is a leafless $\R$-tree with $D$ a subset of $\partial T$ and $G$ a group of isometries of $T$ stabilizing $D$, is a \emph{sparse tuple} if the following two conditions are satisfied.

\begin{itemize}
\item[(S1)]
The group $G$ is free.
\item[(S2)]
If $a \in  \partial T \setminus D$, then the stabilizer of $a$ consists of hyperbolic isometries and the identity only. 
\end{itemize}

A sparse tuple will be called a \emph{2-transitive} if $G$ induces a 2-transitive action on $D$. 


\subsubsection{Sketch of the construction}\label{section:sketch}

Our construction of a counterexample consists of various steps enriching a given sparse tuple $(T,\dd,D,G)$. In particular we will make use of the techniques outlined in Section~\ref{section:treesq}. Note that Conditions (S1) and (S2) behave well under the limit operation of Lemma~\ref{lemma:treesq}, hence we obtain sparse tuples as direct limits of sequences of sparse tuples.

Each of these steps will be such that if they are applied to two sparse tuples $(T_1,\dd_1,D_1,G)$ and $(T_2,\dd_2,D_2,G)$ with a $G$-equivariant bijection $\varphi : D_1 \to D_2$, then for the resulting sparse tuples $(T_1',\dd_1',D_1',G')$ and $(T_2',\dd_2',D_2',G')$ the map $\varphi$ extends to a $G'$-equivariant bijection $\varphi': D_1' \to D_2'$.

In Step A we will show how extend $G$ by adding a single isometry. We then use this in Step B to produce a 2-transitive sparse tuple. 

Step C produces fixed ends in $D$ for those isometries in $G$ for which no non-trivial power fixes more than one end in $D$. Step D repeats step C twice and then constructs a sparse tuple  $(T',\dd',D',G)$ such that $D' = \partial T'$. 

Step E combines Steps B and D to produce a two-transitive sparse tuple $(T',\dd',D',G)$ with $D' = \partial T'$. 

Consider a pair of sparse tuples $(T_1,\dd_1,D_1,G)$ and $(T_2,\dd_2,D_2,G)$ with a $G$-equivariant bijection $\varphi : D_1 \to D_2$ which does not extend to $T_1$ and $T_2$. (An example of this is easily constructed, especially if $G$ is the trivial group.) By Step E one obtains two 2-transitive tuples $(T_1',\dd_1',\partial T_1',G')$ and $(T_2',\dd_2',\partial T_2',G')$ with a $G'$-equivariant bijection between $\partial T_1'$ and $\partial T_2'$ which does not extend to $T_1'$ and $T_2'$, which provides the desired counterexample.

\subsubsection{Step A - Extending the group of isometries}

Let $(T,\dd,D,G)$ be a sparse tuple. Pick three pairwise different ends $a,b$ and $c$  in $D$ such that $b$ and $c$ are in different orbits of the stabilizer $G_a$. We will extend the sparse tuple such there exists an isometry $t$ mapping $b$ to $c$ while fixing the end $a$.

Fix a free generating set $H$ of $G$.  Parametrize the union of $H$ and a symbol $t$ by an index set $J$ and a map $h : J \to H \cup \{t\}$. 

We define the (partial) action of $t$ on $T$ as mapping the apartment $(b,a)$ to $(c,a)$ while fixing the intersection of both pointwise (and the partial action of $t^{-1}$ as mapping $(c,a)$ to $(b,a)$ accordingly).

\paragraph{Construction of the intermediate objects}\label{section:con-int.}

In this section we construct trees $(T_i,\dd_i)$ and  subsets $D_i$ of $\partial T_i$ ($i \in \N$) recursively. These constructions will be such that $(T_{i},\dd_{i})$ is a subtree of $(T_{i+1},\dd_{i+1})$ for $i \in \N$.

We define $(T_0,\dd_0)$ to be the tree $(T,\dd)$ and $D_0$ equal to $D$.

We now construct $(T_{i+1},\dd_{i+1})$ provided that we already did construct $(T_i,\dd_i)$ for a certain $i \in \N$. For every $j \in J$ let $S_i^j$ be the subtree of $T_i$ on which  $h(j)$ is defined and $\widetilde{S}_i^j$ be the subtree on which $h(j)^{-1}$ is defined.

Clearly the subtrees $S_i^j$ and $(S_i^j)^{h(j)}$ of $T_i$ are isometric, hence one can glue $T_i$ to a copy $T_i^j$ of $T_i$  along the subtree $(S_i^j)^{h(j)}$ of $T_i$ and the subtree of $T_i^j$ corresponding to $S_i^j$. Similarly one can glue $T_i$ to a copy $\widetilde{T}_i^j$  along the subtree $(\widetilde{S}_i^j)^{h_j^{-1}}$ of $T_i$ and the subtree of $\widetilde{T}_i^j$ corresponding to $\widetilde{S}_i^j$.

Using the transfinitely recursive method of successively glueing trees together outlined in Section~\ref{section:treesq}, we can repeat this gluing for every $j \in J$. We denote the resulting tree by $(T_{i+1},\dd_{i+1})$. On this tree we can extend the action of $h(j)$ (for each $j \in J$) by mapping $T_i$ to $T_i^j$ via the canonical isometry, and the action of $h(j)^{-1}$ by mapping $T_i$ to $\widetilde{T}_i^j$, again via the canonical isometry.

\paragraph{The limit object and its properties.}

The algorithm from Section~\ref{section:con-int} produces a sequence of trees $(T_0,\dd_0)$, $(T_1,\dd_1)$, $\dots$ to which we can apply Lemma~\ref{lemma:treesq} to combine these into a tree $(T',\dd')$. Each of the $h(j)$ ($j \in J$) act as isometries of this tree, and hence generate a group $G'$ of isometries acting on $T'$. Let $D' = \{fd \vert f \in G', d \in D\}$. 


The claim is that $(T',\dd',D', G')$ is a sparse tuple.

Note that $T$ is a subtree of $T'$ which intersects each orbit of $G'$ on the points of $T'$. We define the \emph{depth} $\delta(x)$ of a point $x \in T'$ as  $\min \{ i \in \N \vert x \in T_i  \}$. We now list some observations on this notion of depth.

\begin{lemma}\label{lemma:step}
If $g \in H \cup \{ t\}$ and $x \in T'$, then $\vert \delta(gx) - \delta(x) \vert \leq 1$. Moreover $\delta(gx) = \delta(x)$ implies $\delta(x) = 0$.
\end{lemma}
\proof
This holds as $T_{i+1}$ extends $T_i$ in such a way that each $g \in H \cup H^{-1} \cup \{ t, t^{-1}\} $ is defined on $T_i$ with $g(T_i) \subset T_{i+1}$.
\qed

\begin{lemma}\label{lemma:valley}
Let $a_n \dots a_1$ ($n \in \N$) be a reduced non-trivial word in the generators $H \cup \{t\}$ and their inverses. Let $g$ be the corresponding isometry of $G'$ and $x$ a point of $T'$. If $\delta(a_1 x)  = \delta(x) + 1$, then $\delta(g x)  = \delta(x) +n$. 
\end{lemma}
\proof
This follows from Lemma~\ref{lemma:step} and the construction of $T'$.
\qed

%

In the next few lemmas we use depth to study the structure of $G'$.

\begin{lemma}\label{lemma:bound}
Let $a_n \dots a_1$ ($n \in \N$) be a cyclically reduced non-trivial word in the generators $H \cup \{t\}$ and their inverses. Let $g$ be the corresponding isometry of $G'$, then $g$ only fixes points of $T$.
\end{lemma}
\proof
Suppose that $g$ fixes a point $x$ of $T'$. The subsequent applications of the $a_i$ ($i\in \{1, \dots, n \}$) to the point $x$, either preserve the depth, increase it by one or decrease it by one at each step by Lemma~\ref{lemma:step}, starting and ending at $\delta(x)$. If this depth increments by one at a certain step, then a cyclic permutation of $a_n \dots a_1$ will fix some point $y$ while at the same time, by Lemma~\ref{lemma:valley}, it increases the depth of $y$ by $n$. Hence the depth stays constant at each step. From Lemma~\ref{lemma:step} it follows that $x \in T$.
\qed

\begin{lemma}\label{lemma:free}
The group $G'$ is a free group with free generating set $H \cup \{t\}$.
\end{lemma}
\proof
If $H \cup \{t\}$ is not a free generating set then there would exist a cyclically reduced non-trivial word $a_n \dots a_1$ ($n \in \N$) in the generators $H \cup \{t\}$ and their inverses which corresponds with the identity isometry of $T'$. However, the corresponding isometry would act non-trivially on points in $T' \setminus T$ by Lemma~\ref{lemma:bound}. This proves the lemma.
\qed

\begin{prop}
The tuple $(T',\dd',D', G')$ is sparse.
\end{prop}
\proof
Condition (S1) is shown in Lemma~\ref{lemma:free}. For Condition (S2) suppose by way of contradiction that there exists a non-trivial elliptic isometry $g \in G'$ which fixes an end $d \in \partial T' \setminus D'$. In particular it will fix a ray with end $d$ pointwise.

Let $a_n \dots a_1$ ($n \in \N$) be a reduced word in the generators $H \cup \{t\}$ and their inverses representing $g$. Without loss of generality we may assume that $a_n \dots a_1$ is cyclically reduced as its cyclic permutations fix ends in the $G'$-orbit of $d$, which is disjoint from $D'$ as the latter is stabilized by $G'$. 

Lemma~\ref{lemma:bound} states that $g$ fixes only elements of $T$, hence $d$ is an end of $T$. Similarly the ends $a_1 d$, $a_2 a_1 d$, \dots lie completely in $\partial T \setminus D$ by considering the cyclic permutations of $a_n \dots a_1$. As $t$ and $t^{-1}$ both map $\partial T$ completely outside of $\partial T$ except for $a$, $b$ and $c$ which all lie in $D$, none of the $a_i$ ($i\in \{1, \dots n\}$) equal $t$ or $t^{-1}$. It follows that $g \in G$,  hence Condition (S2) for $(T,\dd,D,G)$ yields a contraction. We conclude that $(T',\dd',D', G')$ is a sparse tuple. \qed

%
%


\paragraph{Extending pairs of sparse tuples.}\label{section:extendpairs}

Let $(T_1,\dd_1,D_1,G)$ and $(T_2,\dd_2,D_2,G)$ be two sparse tuples together with a $G$-equivariant bijection $\varphi : D_1 \to D_2$. Suppose we apply Step A to $(T_1,\dd_1,D_1,G)$ adding an isometry $t$ mapping the apartment $(a,b )$ to the apartment $(a,c)$ for pairwise distinct ends $a,b$ and $c$ in $D_1$ such that $b$ and $c$ are not in the same $G_a$ orbit. Similarly we apply Step A to $(T_2,\dd_2,D_2,G)$ for the ends $\varphi (a)$, $\varphi( b)$ and $\varphi (c)$. In this way we obtain two new sparse tuples $(T_1',\dd_1',D_1',G')$ and $(T_2',\dd_2',D_2',G')$.

We now consider the action of $G'$ on the ends in $D_1$.

\begin{lemma}\label{lemma:no_dep}
If $d$ and $e$ are two ends of $D_1$, then the maps in $G'$ mapping $d$ to $e$ and those that map $\varphi (d)$ to $\varphi(e)$ are the same. In particular the $G'$-stabilizers of $d$ and $\varphi(d)$ are identical.
\end{lemma}
\proof
Let $g  \in G'$ an element that maps the end $d$ to $e$, represented by a reduced word $a_n \dots a_1$ ($n \in \N$) in the generators $H \cup \{t\}$ and their inverses.

Consider the images of $d$ under $a_1$, $a_2a_1$, $\dots$ up to $g=a_n \dots a_1$. If each of these images lies in $D_1$, then $g$ also maps $\varphi(d)$ to $\varphi(e)$ via the $G$-equivariant bijection $\varphi$ and the equivariance for the partial action of $t$ on $T_1$ and $T_2$.

Now assume that $a_j \dots a_1 d$ is no longer in $D_1$ and that $j \in \N$ is minimal with this property. Then $a_j$ equals either $t$ or $t^{-1}$ and maps some ray with end $d$ completely contained in $T_1$ to a ray disjoint from $T_1$. From Lemma~\ref{lemma:valley} it follows that the $g$ maps this ray to a ray consisting completely of depth non-zero points, hence it cannot have as end $e$. 

This implies that the elements in $G'$ mapping $d$ to $e$ also map $\varphi(d)$ to $\varphi(e)$. Analogously one proves the other inclusion. 
\qed

We can now extend $\varphi$ to $D_1'$ as follows. Let $m$ be an end in $D_1'$.  As each orbit of $G'$ on $D_1'$ intersects $D_1$ (by construction of $D_1'$) there exists a $g \in G'$ mapping an end $d \in D_1$ to $m$.  We then define $\varphi'(m)$ as the end $g \varphi(d) \in D_2'$. The next  two lemmas show that $\varphi'$ is well-defined and a $G'$-equivariant bijection between $D_1'$ and $D_2'$, as is required by Section~\ref{section:sketch}.

\begin{lemma}
The map $\varphi'$ is well-defined.
\end{lemma}
\proof
In order to show that $\varphi'$ is well-defined we have prove that it is independent of the choice of the end $d \in D_1$ and the element $g \in G'$. Suppose $d,e \in D_1$ and $g,h \in G'$ are such that $gd=m$ and $ he = m$ with $m \in D_1'$. We now need that $ h \varphi(e) = g \varphi(d) $.

As $d$ and $e$ have to be in the same $G'$-orbit, there exists a $g' \in G'$ mapping $d$ to $e$. Lemma~\ref{lemma:no_dep} states that $g'$ also maps $\varphi(d)$ to $\varphi(e)$. Because $g^{-1} h g'$ is in the $G'$-stabilizer of $d$, it is also in the $G'$-stabilizer of $\varphi(d)$, again by Lemma~\ref{lemma:no_dep}. So 
\begin{align*}
g^{-1}  h g'  \varphi(d) = \varphi(d)  \Leftrightarrow & h g'  \varphi(d) = g \varphi(d) \\
 \Leftrightarrow& h \varphi(e) = g \varphi(d), 
\end{align*}
which is what we needed.
\qed

\begin{lemma}
The map $\varphi'$ is a $G'$-equivariant bijection between $D_1'$ and $D_2'$.
\end{lemma}
\proof
By switching the roles of $D_1'$ and $D_2'$ we see that $\varphi'$ is a bijection between $D_1'$ and $D_2'$. Let $g \in G'$ and $d \in D_1$ be such that $g d = m $ with $m \in D_1'$. Now let $h \in G'$, in order to have a $G'$-equivariant bijection we need that $ \varphi'(h m) = h \varphi'(m)$.

This is the case as 
$$
h \varphi'(m) = h g \varphi (d)  =  \varphi' ( h g d) = \varphi' (h m)
$$
by the definition of $\varphi'$. \qed

\subsubsection{Step B - 2-transitive sparse tuples}

In Step A we enlarged a sparse tuple $(T,\dd,D,G)$ by adding isometries mapping a certain pair of ends $(a,b)$ to a pair $(a,c)$.

By repeating Step A recursively for some large enough ordinal and applying Lemma~\ref{lemma:treesq} we are able to enlarge the sparse tuple $(T,\dd,D,G)$ to a sparse tuple $(T',\dd',D',G')$ where $G'$ acts two-transitive on the set of ends $D'$. 

The results of Section~\ref{section:extendpairs} hold at each step of the recursion, so they also hold for this Step B.

\subsubsection{Step C - Enriching the fixed ends}
Let $(T,\dd,D,G)$ be a sparse tuple. In this step  we will produce fixed ends in $D$ for those isometries in $G$ for which no non-trivial power fixes more than one end in $D$.

We denote by $L$ the tree consisting of a single ray.

Consider the maximal cyclic subgroups of $G$ for which no non-trivial element fixes more than one end in $D$, and pick a generator of such a subgroup in each conjugacy class of these. Denote the set of the elements picked this way by $H$. Note that every cyclic group of $G$ infinite and is contained in a unique maximal one as the group $G$ is free.

%

We partition this set into two subsets $H_e$ and $H_h$, consisting of respectively the elliptic and hyperbolic isometries in $H$. 

We first discuss the set $H_e$. For every $g \in H_{e}$ pick a fixed point  $x_g$ of $g$ in $T$. For every left coset $V$ in $ G/  \langle g \rangle $ glue a copy $L_V$ of $L$ to $T$ along the point $V x_g$ (such that the origin of $L_V$ is identified with $Vx_g$). Note that $Vx_g$ is a unique point as $x_g$ is fixed by $g$ and hence also by $\langle g \rangle$. Doing this recursively for each $g \in H_e$ we obtain a tree $(T',\dd')$. 

The action of $G$ on $T$ can be extended to $T'$ by letting $h \in G$ map $L_V$ (with $V \in G / \langle g \rangle$, $g \in H_e$) to $L_{hV}$ via the canonical isometry. Let $D'$ be the union of $D$ and the ends $d_g$ of the glued $L_V$'s (again for $V \in G / \langle g \rangle$ and $g \in H_e$). As $\partial T \setminus D = \partial T' \setminus D'$ one readily observes that $(T',\dd',D',G)$ is again a sparse tuple.

We now discuss $H_h$. For every $g \in H_h$ pick an end $d_g \in \partial T' \setminus D'$ fixed by $g$. Such an end always exists by Lemma~\ref{lemma:dich} and the construction of $H_h$. We now extend $D'$ by adding the orbit of $d_g$ under $G$ for every $g \in H_h$ to it, resulting in a subset $D'' \subset \partial T'$. Clearly $(T',\dd',D'', G)$ is again a sparse tuple

Important to note is that for every $g \in H$, elliptic or hyperbolic, the stabilizer of the end $d_g$ is exactly $\langle g \rangle$. In the elliptic case this follows from the construction, in the hyperbolic case from (S2) and maximality of the cyclic subgroups. So if $(T_1,\dd_1,D_1,G)$ and $(T_2,\dd_2,D_2,G)$ are two sparse tuples with a $G$-equivariant  bijection $\varphi : D_1 \to D_2$, then for the resulting sparse tuples $(T_1',\dd_1',D_1'',G)$ and $(T_2',\dd_2',D_2'',G)$ the map $\varphi$ extends to a $G$-equivariant bijection $\varphi': D_1'' \to D_2''$.

\subsubsection{Step D - Completing the set of ends}
The first part of this step is to apply Step C twice to the given sparse tuple $(T,\dd,D,G)$ to obtain a sparse tuple $(T',\dd',D',G)$. This ensures that every isometry in $G$ has a non-trivial power fixing at least two ends in $D$. The tuple $(T',\dd',\partial T',G)$, where we extend $D'$ maximally, is sparse for trivial reasons. 

Note that if $ d \in \partial T' \setminus D'$, then the stabilizer of $d$ in $G$ is trivial. Hence one can easily extend $G$-equivariant bijections for pairs of sparse tuples.

\subsubsection{Step E - Combining the steps}

We start with a sparse tuple $(T,\dd, D, G)$. On this sparse tuple we apply Steps $D$ and $B$ alternately, repeating this transfinitely recursively (using Section~\ref{section:treesq}) for the first uncountable ordinal $\omega_1$. As limit object we obtain a sparse tuple $(T',\dd',D',G')$.

As the cofinality of $\omega_1$ equals its own cardinality $\aleph_1$, every countable subset of $\omega_1$ (considered as a well-ordered set) is bounded from above by some element of $\omega_1$.  In particular this yields that $\partial T' = D'$ (as a ray, so also its end, can be defined by as the convex closure of a countable subset of this ray together with Step D), and that $G'$ acts two-transitively on $D'$ (from Step B).

One has that $G$-equivariant bijections extend well, as this is the case for Steps B and D.

%

\end{document}